\documentclass[12pt]{article}
\usepackage{graphicx}
\usepackage{subfig}
\usepackage[a4paper, total={7in, 10.5in}]{geometry}
\usepackage{xcolor}

\usepackage{amsmath}
\usepackage{amsfonts}
\usepackage{amssymb}
\usepackage{amsthm}
\usepackage{hyperref}

\usepackage[ruled, linesnumbered]{algorithm2e}

\usepackage{svg}

\hypersetup{hidelinks}

\usepackage{listings}
\usepackage{doi}

\newcommand{\NN}{\mathbb N}
\newcommand{\PP}{\mathbb P}

\newcommand{\ZZ}{\mathbb Z}

\theoremstyle{plain}
\newtheorem{theorem}{Theorem}[section]

\theoremstyle{definition}

\theoremstyle{remark}
\newtheorem{remark}[theorem]{Remark}

\title{Computing all minimal Markov bases in Macaulay2}
\author{Oliver Clarke, Alexander Milner}
\date{February 2025}

\begin{document}

\maketitle

\begin{abstract}
    We introduce the package \textit{allMarkovBases} for \textit{Macaulay2}, which is used to compute all minimal Markov bases of a given toric ideal. The package builds on functionality of \textit{4ti2} by producing the \textit{fiber graph} of the toric ideal. The package uses this graph to compute other properties of the toric ideal such as its indispensable set of binomials as well as its universal Markov basis.
\end{abstract}

\section{Introduction}

Let $k$ be a field and fix a matrix $A \in \ZZ^{d \times n}$ such that $\ker(A) \cap \NN^n = \{0\}$, which we call a configuration matrix. The toric ideal $I_A \subseteq k[x_1, \dots, x_n]$ is the ideal generated by the binomials $x^u - x^v$ where $u,v \in \NN^n$ such that $u - v \in \ker(A)$. Note that here we use multi-index notation, so $x^u$ denotes the monomial $x_1^{u_1}x_2^{u_2} \cdots x_n^{u_n}$.
For an introduction to toric ideals, see \cite{sturmfels1996grobner}.
A \textit{(minimal) Markov basis} of $A$ is an (inclusion-minimal) set of generators $M$ of $I_A$ and we identify $M$ with its set of binomial exponents $\{u-v : x^u - x^v \in M\}$. The seminal paper of Diaconis and Sturmfels \cite{DS_1998} provides a connection between Markov bases of toric ideals and sampling from conditional distributions. The research landscape of Algebraic Statistics has since seen a great deal of development, particularly in the study of Markov bases. For an introduction to Markov bases in Statistics we refer to \cite{Kos_2020} and \cite{ALP_2024}, and for a general introduction to Algebraic Statistics we refer to \cite{sullivant2023algebraic}.

The current version of \textit{Macaulay2} \cite{M2} is shipped with \textit{4ti2} \cite{4ti2}, which is a powerful library for performing computations with toric ideals, polyhedra, and lattices. Users of \textit{Macaulay2} may be familiar with its interfacing package \textit{FourTiTwo}, which provides the function \texttt{toricMarkov} that computes a single Markov basis from a given matrix.
However, this package does not provide a way to study all possible minimal Markov bases of a given matrix $A$. So we introduce the package \textit{allMarkovBases}, which comes with a function \texttt{markovBases} that returns all minimal Markov bases. We do this by computing the \textit{fiber graphs} of $A$, which is a sequence of graphs whose connectivity can be used to determine the Markov bases of $A$. See Section~\ref{sec: fiber graph}. Our package exploits the fiber graph to compute important families of binomials of the toric ideal, which cannot be computed directly with \textit{FourTiTwo}. This includes the \textit{indispensable set} $S(A)$ and \textit{universal Markov basis} $U(A)$, which are, respectively, the intersection and union of all minimal Markov bases. These families of elements have alternative combinatorial descriptions and prominently feature in the study of Markov bases. See, for example, \cite{
Charalambous2014Markov, 
Kosta2023Strongly, 
clarke2024distance,
CTV_2017}. If the number of Markov bases is very large, then our package allows for the uniform sampling from the set of Markov bases.

\subsection{Outline with examples}

The matrix $A = \begin{pmatrix}
    7 & 8 & 9 & 10
\end{pmatrix}$ has four minimal Markov bases, which are computed as follows.
\begin{lstlisting}
i1 : needsPackage "allMarkovBases";

i2 : A = matrix "7,8,9,10";

i3 : countMarkov A
o3 = 4

i4 : markovBases A
o4 = {| -1 2  -1 0  |, | -1 2  -1 0  |, | -1 2  -1 0  |, | -1 2  -1 0  |}
      | -1 1  1  -1 |  | -1 1  1  -1 |  | -1 1  1  -1 |  | -1 1  1  -1 |
      | 0  -1 2  -1 |  | 0  -1 2  -1 |  | 0  -1 2  -1 |  | 0  -1 2  -1 |
      | 4  0  -2 -1 |  | 4  0  -2 -1 |  | 4  -1 0  -2 |  | 4  -1 0  -2 |
      | 3  1  -1 -2 |  | 3  1  -1 -2 |  | 3  1  -1 -2 |  | 3  1  -1 -2 |
      | 3  0  1  -3 |  | 2  2  0  -3 |  | 3  0  1  -3 |  | 2  2  0  -3 |
\end{lstlisting}

On the other hand, the matrix $A' = \begin{pmatrix}
    51 & 52 & 53 & 54 & 55 & 56
\end{pmatrix}$ has $24300$ minimal Markov bases. It may not be feasible to work with all of them so we can produce random samples as follows.
\begin{lstlisting}
i5 : A' = matrix "51,52,53,54,55,56";

i6 : countMarkov A'
o6 = 24300

i7 : randomMarkov A'
o7 = | 8  4  0  0  0  -11 |
     | -1 2  -1 0  0  0   |
     | -1 1  1  -1 0  0   |
     | -1 0  2  0  -1 0   |
     | -1 1  0  1  -1 0   |
     | 0  -1 1  1  -1 0   |
     | -1 1  0  0  1  -1  |
     | 0  0  -1 2  -1 0   |
     | 0  -1 0  2  0  -1  |
     | 0  0  -1 1  1  -1  |
     | 0  0  0  -1 2  -1  |
     | 12 0  -1 0  -1 -9  |
     | 11 1  0  -1 -1 -9  |
     | 11 0  1  -1 0  -10 |
     | 9  3  0  0  -1 -10 |
\end{lstlisting}

The \textit{indispensable set} $S(A')$ and the \textit{universal Markov basis} $U(A')$ are computed as follows. We have abridged the output of \texttt{toricUniversalMarkov A'} to save space. We note that the indispensable set and universal Markov bases are read from the fiber graph, so the methods do not require the computation of all 24300 minimal Markov bases.

\begin{lstlisting}
i8 : toricIndispensableSet A'
o8 = | -1 2 -1 0  0 0  |
     | -1 1 1  -1 0 0  |
     | 0  0 -1 1  1 -1 |
     | 0  0 0  -1 2 -1 |

i9 : toricUniversalMarkov A'
o9 = | 10 0  2  0  0  -11 |
     | 8  4  0  0  0  -11 |
      ...
     | 9  3  0  0  -1 -10 |
     | 11 0  0  1  -1 -10 |
              33       6
o9 : Matrix ZZ   <-- ZZ
\end{lstlisting}

An example from Algebraic Statistics is the Segre embedding of a product of projective spaces: $\PP^2 \times \PP^2 \times \PP^2$, whose toric ideal is given by the matrix
\[
A'' = \begin{bmatrix}
    1&1&1&1&1&1&1&1&1&1&1&1&1&1&1&1&1&1&1&1&1&1&1&1&1&1&1\\
    1&1&1&1&1&1&1&1&1&0&0&0&0&0&0&0&0&0&0&0&0&0&0&0&0&0&0\\
    0&0&0&0&0&0&0&0&0&1&1&1&1&1&1&1&1&1&0&0&0&0&0&0&0&0&0\\
    1&1&1&0&0&0&0&0&0&1&1&1&0&0&0&0&0&0&1&1&1&0&0&0&0&0&0\\
    0&0&0&1&1&1&0&0&0&0&0&0&1&1&1&0&0&0&0&0&0&1&1&1&0&0&0\\
    1&0&0&1&0&0&1&0&0&1&0&0&1&0&0&1&0&0&1&0&0&1&0&0&1&0&0\\
    0&1&0&0&1&0&0&1&0&0&1&0&0&1&0&0&1&0&0&1&0&0&1&0&0&1&0  
\end{bmatrix}.
\]
In this case, we were able to compute the number of Markov bases, the indispensable set, and universal Markov basis in under a second. 
\begin{lstlisting}
i10 : A'' = matrix {{1,1,1,1,1,1,1,1,1,1,1,1,1,1,1,1,1,1,1,1,1,1,1,1,1,1,1},
      		{1,1,1,1,1,1,1,1,1,0,0,0,0,0,0,0,0,0,0,0,0,0,0,0,0,0,0},
      		{0,0,0,0,0,0,0,0,0,1,1,1,1,1,1,1,1,1,0,0,0,0,0,0,0,0,0},
      		{1,1,1,0,0,0,0,0,0,1,1,1,0,0,0,0,0,0,1,1,1,0,0,0,0,0,0},
      		{0,0,0,1,1,1,0,0,0,0,0,0,1,1,1,0,0,0,0,0,0,1,1,1,0,0,0},
      		{1,0,0,1,0,0,1,0,0,1,0,0,1,0,0,1,0,0,1,0,0,1,0,0,1,0,0},
      		{0,1,0,0,1,0,0,1,0,0,1,0,0,1,0,0,1,0,0,1,0,0,1,0,0,1,0}};

i11 : elapsedTime countMarkov A''
 -- .306693s elapsed
o11 = 324518553658426726783156020576256

i12 : elapsedTime toricIndispensableSet A'';
 -- .00462137s elapsed
               81       27
o12 : Matrix ZZ   <-- ZZ

i13 : elapsedTime toricUniversalMarkov A'';
 -- .0104826s elapsed
               243       27
o13 : Matrix ZZ    <-- ZZ
\end{lstlisting}
In the above, most of the time is spent computing the fiber graph. Once the graph has been computed, it is stored to avoid re-computation. The above example is particularly amenable to our computations because the fiber graphs are all very small.

\section{How it works}

Throughout this section, we take $A$ to be a $d \times n$ configuration matrix. We explain how our functions work with the simple running example $A = \begin{pmatrix}
    1 & 2 & 3
\end{pmatrix} \in \ZZ^{1 \times 3}$. First, we require a distinguished Markov basis $M$, which is computed by \textit{FourTiTwo} as follows.
\begin{lstlisting}
i1 : needsPackage "allMarkovBases";
i2 : A = matrix "1,2,3";

i3 : M = toricMarkov A
o3 = | 2 -1 0  |
     | 3 0  -1 |
\end{lstlisting}
So the toric ideal $I_A \subseteq k[x,y,z]$ is minimally generated by the binomials $x^2 - y$ and $x^3 - z$. The \textit{affine semigroup} of $A$ is $\NN A := \{Ax \in \ZZ^d : x \in \NN^n\}$ and, for each $t \in \NN A$, the  \textit{$t$-fiber} of $A$ is the subset $\mathcal F_t := \{u \in \NN^n : Au = t\}$. For each binomial $x^u - x^v \in I_A$, we have that $u, v \in \mathcal F_t$ belong to the same fiber of $A$, and we define its \textit{$A$-degree} $\deg_A(x^u - x^v) = t$.
In our running example we have 
\[
\deg_A(x^2 - y) = 2 
\quad \text{and} \quad
\deg_A(x^3 - z) = 3.
\]

By \cite[Theorems~2.6 and 2.7]{CKT_2007} the multiset of $A$-degrees is the same for any minimal Markov basis. For instance, the binomials $x^2 - y,\, xy - z$ form another minimal Markov basis for $I_A$ and and their $A$-degrees are $2$ and $3$ respectively. We say that $t \in \NN A$ is a \textit{generating fiber} for $A$ if there exists an element $x^u - x^v$ of a minimal Markov basis with $A$-degree $\deg_A(x^u - x^v) = t$. The aforementioned results provide a natural procedure to generate all minimal Markov bases using the fiber graph of $A$.

\subsection{Fiber graphs}\label{sec: fiber graph}
Let $t \in \NN^d$. The \textit{fiber graph} $G_t$ has vertex set $\mathcal F_t$ and edge set 
\[
E(G_t) = \{uv : \text{there exists } i \in \{1, \dots, n\} \text{ such that } u_i > 0 \text{ and } v_i > 0\}.
\]
Using our package, the list of fiber graphs of the generating fibers of $A$ are computed as follows.

\begin{lstlisting}
i4 : fiberGraph A
o4 = {Graph{{0, 1, 0} => {}}, Graph{{0, 0, 1} => {}         }}
            {2, 0, 0} => {}         {1, 1, 0} => {{3, 0, 0}}
                                    {3, 0, 0} => {{1, 1, 0}}
\end{lstlisting}

\begin{figure}
    \centering
    \subfloat[\label{fig:2fiber} \centering $A$-degree 2]{{\includegraphics[scale = 0.23]{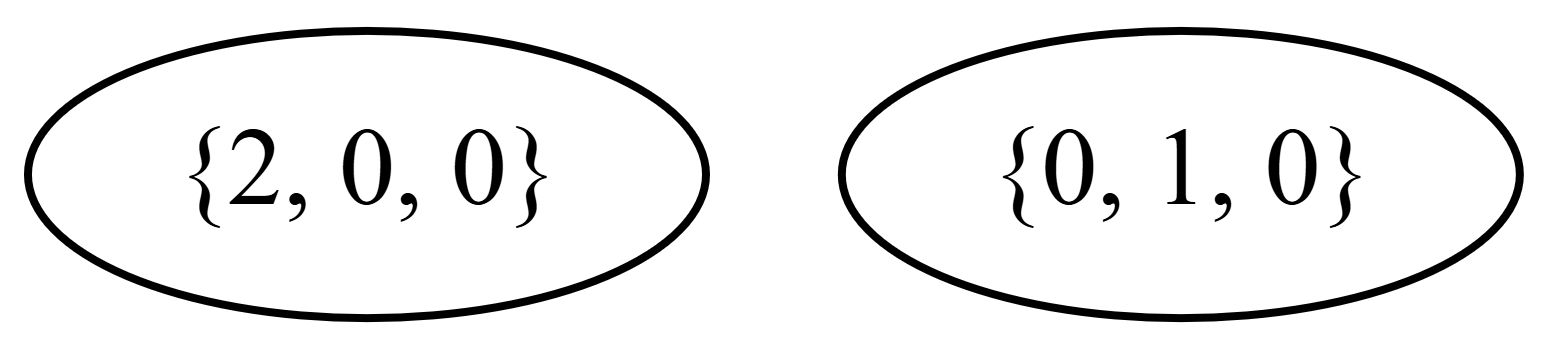} }}
    \qquad
    \subfloat[\label{fig:3fiber} \centering $A$-degree 3]{{\includegraphics[scale = 0.23]{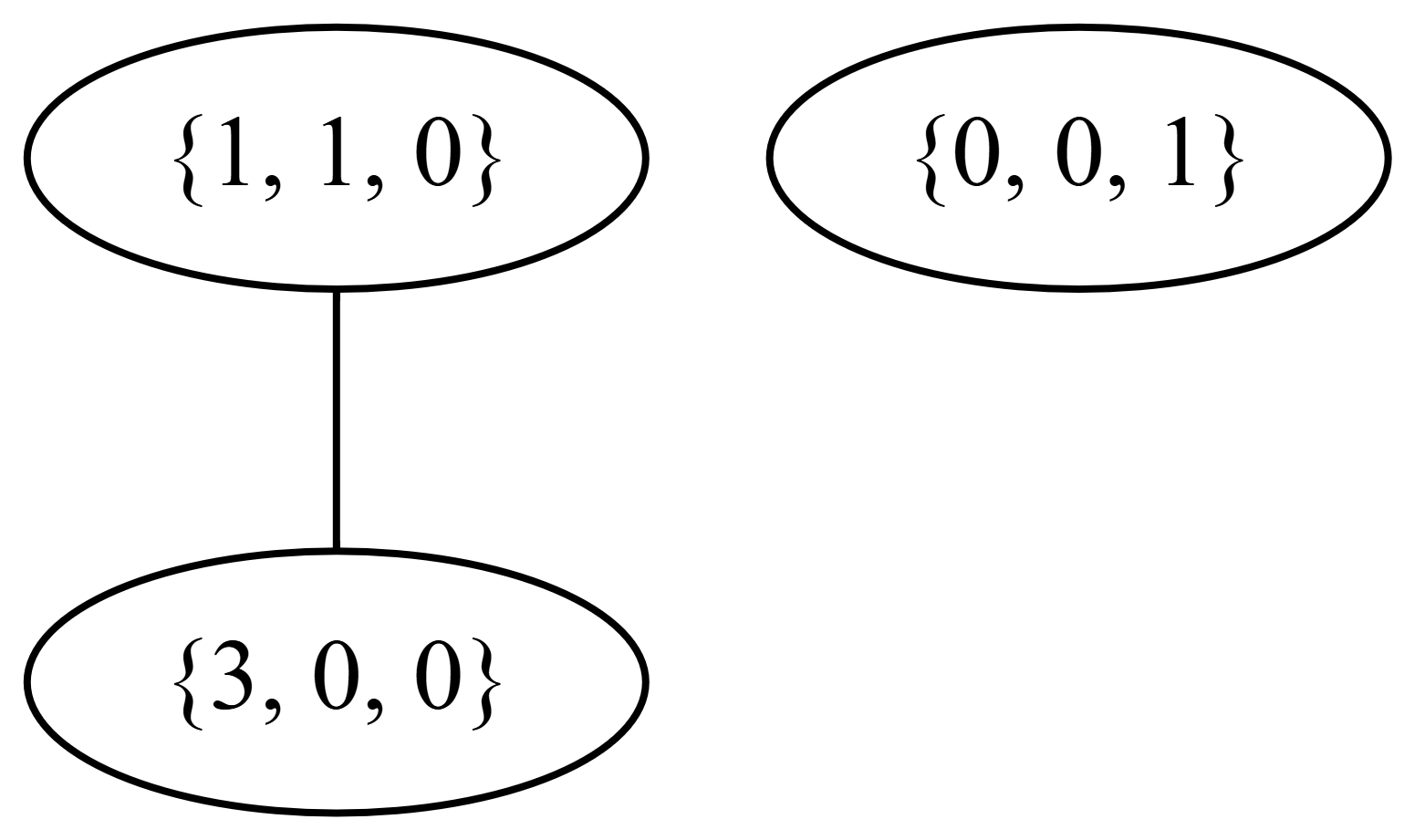} }}
    \caption{Fiber graphs of $A$ for the generating fibers}
    \label{fig:fibers}
\end{figure}

A visualisation of these fibers, produced by the package \textit{Graphs}, is shown in Figure~\ref{fig:fibers}.  Each element $z \in \ker(A)$ is thought of as a \textit{move} that connects pairs of points $u, v \in \NN^n$ whenever $u = v + z$. In our example, the move $(2, -1, 0)$ connects $(2,0,0)$ and $(0,1,0)$ and connects $(1,1,0)$ and $(3,0,0)$.
The \textit{fundamental theorem of Markov bases} \cite[Theorem~3.1]{DS_1998} tells us that the set of Markov bases are in one-to-one correspondence with collections of elements of $\ker(A)$ that \textit{connect} each fiber of $A$. Here, a subset $M \subseteq \ker(A)$ is said to \textit{connect} the $t$-fiber, for some $t \in \NN A$, if the graph $G_{M, t}$ on $\mathcal F_t$ with edges $\{uv : u-v \in M \}$ is connected. Similarly, by \cite[Theorems~2.6 and 2.7]{CKT_2007}, the minimal Markov bases are in one-to-one correspondence with sets of moves that \textit{minimally connect $G_t$} for each generating fiber $t$ of $A$. 

\begin{remark}
    In the literature, the graphs $G_{M,t}$ are often called fiber graphs \cite{Hemmecke2009Computing}. To avoid confusion, for us the term fiber graph will always refer to the graph $G_t$ for some $t$. 
\end{remark}

In our example, we observe that the $2$-fiber in Figure~\ref{fig:2fiber} is connected only by the element $(2,-1,0)$, up to sign. It follows that the corresponding binomial $x^2-y$ is indispensable as there is no other way to connect the $2$-fiber, thus every Markov basis must contain it. On the other hand, the $3$-fiber in Figure~\ref{fig:3fiber} is minimally connected with either the move $(1,1,-1)$ or $(3,0,-1)$. These form the two minimal Markov bases of $A$. This process is carried out by the function \texttt{markovBases}:

\begin{lstlisting}
i5 : markovBases A
o5 = {| 2 -1 0  |, | 2 -1 0  |}
      | 3 0  -1 |  | 1 1  -1 |
\end{lstlisting}

\subsection{The all-Markov algorithm}
Our package computes all Markov bases with an implementation of Algorithm~\ref{alg:allMarkov}. By \cite[Theorems~2.6 and 2.7]{CKT_2007}, this algorithm is correct and returns all minimal Markov bases of $A$.

\medskip

\begin{algorithm}[H]\label{alg:allMarkov}
\SetAlgoLined
\SetNlSty{textit}{}{:}
\SetKwInOut{input}{Input}\SetKwInOut{output}{Output}

\input{$A \in \ZZ^{d \times n}$ a configuration matrix}
\output{$\mathcal M$ the set of minimal Markov bases of $A$}
\BlankLine

Initialise $\mathcal M = \emptyset$ \;

\For{$t \in \NN A$ a generating fiber of $A$}{

Compute $C_t := $ the set of connected components of the fiber graph $G_t$ \;

}

\BlankLine

$\mathcal T := \prod_{t} \{T : \text{spanning tree on vertex set } C_t\}$, where the product is taken over all generating fibers $t \in \NN A$ of $A$

\BlankLine

\For{$\{T_t : \text{spanning tree on } C_t\} \in \mathcal T$}{

\For{each collection $\{f_{T, e} : T \in \mathcal T,\, e \in E(T)\}$ where $f_{T, e}$ is a choice function on $e$}{

$M := \{f_{T,e}(u) - f_{T,e}(v) : T \in \mathcal T,\, uv \in E(T)\}$ Markov basis of $A$\;

$\mathcal M \leftarrow \mathcal M \cup \{M\}$\;
}
}

\caption{All minimal Markov bases of a matrix}
\end{algorithm}

\medskip

Let us give some notes on how the algorithm works. For each fiber graph $G_t$, with $t \in \NN A$ a generating fiber of $A$, we treat its connected components as the vertices of a new graph $\widehat G_t$ for which we choose a spanning tree $T$. Then the union over all generating fibers of the binomials $x^u-x^v$, where $u,\,v$ are in distinct connected components of $G_t$ which are connected by an edge of $T$, forms a minimal Markov basis. Concretely, on line~7, the edge $e \in E(T)$ is a pair of connected components of $C_t$, hence $e$ is set of sets of elements of $\NN^n$. Suppose that $uv \in E(T)$ is an edge of some $T \in \mathcal T$, then the choice function $f_{T,e}$ produces the binomial $x^u - x^v$, which appears in the Markov basis.

\begin{figure}
    \centering
    \includegraphics[scale=0.7]{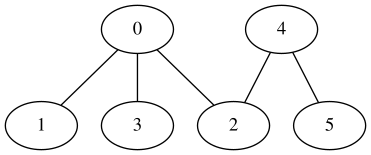}
    \caption{Spanning tree on $\{0,1,\dots,5\}$ associated to the Pr\"ufer sequence $\{0,0,2,4\}$}
    \label{fig:examplePrufer}
\end{figure}

The spanning trees are enumerated using the well-known bijection of Pr\"ufer \cite{Pru_1918}. This bijection between labelled spanning trees on $n$ vertices and \textit{Pr\"ufer sequences} (sequences in $\{0, \dots ,n-1\}$ of length $n-2$) is completely constructive. We have implemented a helpful function \texttt{pruferSequence} that takes a Pr\"ufer sequence and returns the edge set of the corresponding spanning tree. For example, the spanning tree associated to the sequence $\{0,1,2,3\}$ is given computed as follows, see Figure~\ref{fig:examplePrufer}.
\begin{lstlisting}
i6 : pruferSequence {0,0,2,4}
o6 = {set {0, 1}, set {0, 3}, set {0, 2}, set {4, 2}, set {4, 5}}
\end{lstlisting}

By knowing only the sizes of connected components of each fiber graph, it is possible to calculate the number of minimal Markov bases.
Let $t \in \NN A$ be a generating fiber of $A$. 
Let us assume $G_t$ has $n_t$ connected components of size $m_{t,1},m_{t,2}, \dots, m_{t,n_t}$ for some positive integers $m_{t,i}$ with $i \in \{1,2, \dots, n_t\}$. Then, by a slight modification of Pr\"ufer's bijection, see \cite[Theorem~2.9]{CKT_2007}, the number of minimal Markov bases of $A$ is given by \[
\prod_{t} m_{t,1} \cdot m_{t,2} \cdots m_{t,n_t} \cdot (m_{t,1} + m_{t,2} + \dots + m_{t,n_t})^{n_t - 2}.
\]

In our package, we have implemented this as the function \texttt{countMarkov}, which computes the number of minimal Markov bases without enumerating all of them.

\subsection{Variations on the algorithm}

Our package includes variations of Algorithm~\ref{alg:allMarkov}, which allows us to sample from the set of Markov bases and compute the indispensable set and universal Markov basis.

\medskip
\noindent \textbf{Random sampling.} The function \texttt{randomMarkov} produces a random Markov basis as follows. In lines~6 and 7 of Algorithm~\ref{alg:allMarkov}, we sample uniformly an element of $\mathcal T$ and a random collection of choice functions respectively. This produces a single uniformly distributed minimal Markov basis of $A$. The optional argument \texttt{NumberOfBases} specifies how many samples to produces.

\begin{lstlisting}
i7 : randomMarkov(A, NumberOfBases => 1)
o7 = | 2 -1 0  |
     | 3 0  -1 |
\end{lstlisting}

\medskip
\noindent \textbf{Indispensable set and universal Markov basis.}
The universal Markov basis $U(A)$ of the matrix $A$ is the union of its minimal Markov bases. Every minimal Markov basis arises from Algorithm~\ref{alg:allMarkov}, so we may compute $U(A)$ with a slight modification. Specifically, our package provides the function \texttt{toricUniversalMarkov}, which implements the following modification to the algorithm.

\medskip

\begin{algorithm}[H]\label{alg:universalMarkov}
\SetAlgoLined
\SetNlSty{textit}{}{:}
\SetKwInOut{input}{Input}\SetKwInOut{output}{Output}

\input{$A \in \ZZ^{d \times n}$ a configuration matrix}
\output{$U(A)$ the universal Markov basis of $A$}
\BlankLine

Initialise $U(A) = \emptyset$ \;

\For{$t \in \NN A$ a generating fiber of $A$}{

Compute $C_t := $ the set of connected components of the fiber graph $G_t$ \;

$K_t := $ complete graph on vertex set $C_t$ \;

}

\BlankLine

$K := \bigsqcup_t K_t,$ where the union is taken over all generating fibers $t \in \NN A$ of $A$

\BlankLine

\For{$\{u, v\} \in E(K)$}{

\For{each choice function $f_{uv}$ on $\{u, v\}$}{

$U(A) \leftarrow U(A) \cup \{f_{uv}(u) - f_{uv}(v)\}$ \;

}
}

\caption{Universal Markov basis of a matrix}
\end{algorithm}

By \cite[Corollary
2.10]{CKT_2007}, the indispensable elements can be detected directly from the fiber graphs. The function \texttt{toricIndispensableSet} determines the indispensable elements by finding each fiber graph $G_t$ with exactly two vertices that are each singletons. In our running example, the universal Markov basis and indispensable set are computed as follows.

\begin{lstlisting}
i8 : toricUniversalMarkov A
o8 = | 2 -1 0  |
     | 3 0  -1 |
     | 1 1  -1 |

i9 : toricIndispensableSet A
o9 = | 2 -1 0 |
\end{lstlisting}

\medskip
\noindent \textbf{Polynomial output.}
For each method: \texttt{markovBases}, \texttt{randomMarkov}, \texttt{toricUniversalMarkov}, \texttt{toricIndispensableSet},  we may specify a polynomial ring. The function returns an ideal whose generating set, accessible with \texttt{gens}, is the corresponding set of binomials.

\begin{lstlisting}
i10 : R = QQ[x, y, z];

i11 : markovBases(A, R)
               2       3               2
o11 = {ideal (x  - y, x  - z), ideal (x  - y, x*y - z)}

i12 : toricIndispensableSet(A, R)
             2
o12 = ideal(x  - y)

i13 : gens toricUniversalMarkov(A, R)
o13 = | x2-y x3-z xy-z |
\end{lstlisting}

% \printbibliography[
% title={Bibliography},
% heading=bibintoc
% ]

\bibliographystyle{alpha}
\bibliography{bib}

\end{document}